\documentclass[leqno,12pt]{article}

\usepackage{latexsym}
 \usepackage{graphicx}
\usepackage{amsmath}
\usepackage{amssymb}
\usepackage{mathtools}
\usepackage{algorithm}
\usepackage{algorithmic}
\usepackage{cite}
\usepackage{color}
\usepackage{subfig}

\usepackage[margin=1.3in, top=1.3in, bottom=1.3in]{geometry}

\newcommand{\nc}{\newcommand}
\nc{\nt}{\newtheorem}
\nt{thm}{Theorem}[section]
\nt{cor}[thm]{Corollary}
\nt{prop}[thm]{Proposition}
\nt{lem}[thm]{Lemma}
\nt{defn}[thm]{Definition}
\nt{rem}[thm]{Remark}
\nt{exa}[thm]{Example}
\nt{ass}[thm]{Assumption}
\nt{alg}[thm]{Algorithm}
\nt{obs}[thm]{Observation}
\nt{que}[thm]{Question}
\nc{\ip}[2]{\mbox{$\langle #1,#2 \rangle$}}
\nc{\pf}{\noindent{\bf Proof\ \ }}
\nc{\finpf}{\hfill{$\Box$}\linespace}
\nc{\linespace}{\vspace{\baselineskip} \noindent}
\nc{\R}{{\bf R}}
\nc{\D}{{\bf D}}
\nc{\T}{{\bf T}}
\nc{\C}{{\bf C}}
\nc{\M}{{\mathcal M}}
\nc{\Rn}{{\bf R}^n}
\nc{\Cn}{{\bf C}^n}
\nc{\Hn}{{\bf H}^n}
\nc{\Mn}{{\bf M}^n}
\nc{\bx}{\bar{x}}
\nc{\by}{\bar{y}}
\nc{\inT}{\mbox{\rm int}\,}
\nc{\cl}{\mbox{\rm cl}\,}
\nc{\gph}{\mbox{\rm gph}\,}
\nc{\argmin}{\mbox{\rm argmin}\,}

\def\tto{\;{\lower 1pt \hbox{$\rightarrow$}}\kern -12pt
           \hbox{\raise 2.8pt \hbox{$\rightarrow$}}\;}
\newenvironment{myequation}{\setcounter{equation}{\value{thm}}
   \begin{equation}}{\addtocounter{thm}{1}\end{equation}}

\nc{\bmye}{\begin{myequation}}
\nc{\emye}{\end{myequation}}

\begin{document}
\title{
Disk matrices and the proximal mapping \\
for the numerical radius
}
\author{
X.Y. Han
\thanks{ORIE, Cornell University, Ithaca, NY 14853, USA.}
\and
A.S. Lewis
\thanks{ORIE, Cornell University, Ithaca, NY 14853, USA.
\texttt{people.orie.cornell.edu/aslewis} \hspace{2mm} \mbox{}
Research supported in part by National Science Foundation Grant DMS-1613996.}
}
\date{\today}
\maketitle

\begin{abstract}
Optimal matrices for problems involving the matrix numerical radius often have fields of values that are disks, a phenomenon associated with partial smoothness.  Such matrices are highly structured:  we experiment in particular with the proximal mapping for the radius, which often maps $n$-by-$n$ random matrix inputs into a particular manifold of disk matrices that has real codimension $2n$.  The outputs, computed via semidefinite programming, also satisfy an unusual rank property at optimality.   
\end{abstract}
\medskip

\medskip


\section{Introduction}
For any matrix $A$ in the vector space $\Mn$ of $n$-by-$n$ complex matrices (for $n > 1$), the field of values (or numerical range) is the set
\[
       W(A) ~=~ \{u^*Au : u\in\Cn, \|u\|_2 =1 \}.
\]
It is a compact convex set, by the Toeplitz-Hausdorff Theorem.  The numerical radius of $A$ is 
\[
      r(A) ~=~ \max \{ |z| : z \in W(A) \}.
\]
The numerical radius is a norm on $\Mn$, so in particular vanishes only for the zero matrix.
As a consequence of the ``power inequality''
\[
r(A^k) \le \big(r(A)\big)^k \qquad (k=0,1,2,\ldots),
\]
it provides a simple measure of the transient stability of the discrete-time dynamical system $x_{k+1} = Ax_k$.  For a broad discussion, see \cite[\S 25]{HandbookLA}.

Given its role as a measure of system stability, we may seek to optimize the numerical radius over some given set of matrices.  Experiments reported in \cite{lewis-overton-disks} suggest empirically that the optimal solution of a random feedback-control-type problem is often a {\em disk matrix\/}: its field of values is a disk centered at zero.  At first sight, this observation is surprising:  disk matrices have been studied extensively \cite{tam-yang}, and although random matrices in $\Mn$ often have field of values that approximate disks, especially for large $n$ (see \cite{collins-random}), the true disk-matrix property is quite restrictive.  In particular, any disk matrix must be singular, and in fact \cite{Cro16} its zero eigenvalue must have algebraic multiplicity at least two.  

As proved in \cite[Theorem 7.4]{lewis-overton-disks}, however, around disk matrices belonging to a special class (which we call {\em strongly certified}), the numerical radius is {\em partly smooth}.  Specifically, the numerical radius is smooth on the manifold ${\mathcal M}$ consisting of all nearby disk matrices, but along paths crossing ${\mathcal M}$ transversally, its rate of change  increases in a jump discontinuity at ${\mathcal M}$.  For any minimization problem, partial smoothness of the objective relative to some manifold ensures that a random perturbation produces a distribution of random outputs with local support the manifold.  As pointed out in \cite{lewis-overton-disks}, this explains why disk matrices tend to appear as solutions of appropriately randomized optimization problems.  

In the current work, we explore the prevalence of disk matrices in numerical radius optimization first noted in \cite{lewis-overton-disks}.  We make three numerical observations and explore the supporting theory.  We first consider a canonical optimization problem, computing the proximal map for the numerical radius at random matrix inputs and observing that the outputs are often disk matrices.  Secondly, although strong certification appears a stringent requirement, we observe that our output disk matrices are indeed always strongly certified.  We therefore conjecture a converse of \cite[Theorem 7.4]{lewis-overton-disks}:  that the numerical radius is partly smooth around a manifold of disk matrices if and only if those matrices are strongly certified.

Both the proximal map for the numerical radius and the numerical radius itself can be calculated via semidefinite programming.  We observe thirdly that the disk matrix property corresponds to a rank condition in the semidefinite program for the numerical radius.  We relate this property to strong certification.

\section{A numerical experiment}
\begin{figure}
\centering
\includegraphics[width=1.0\textwidth]{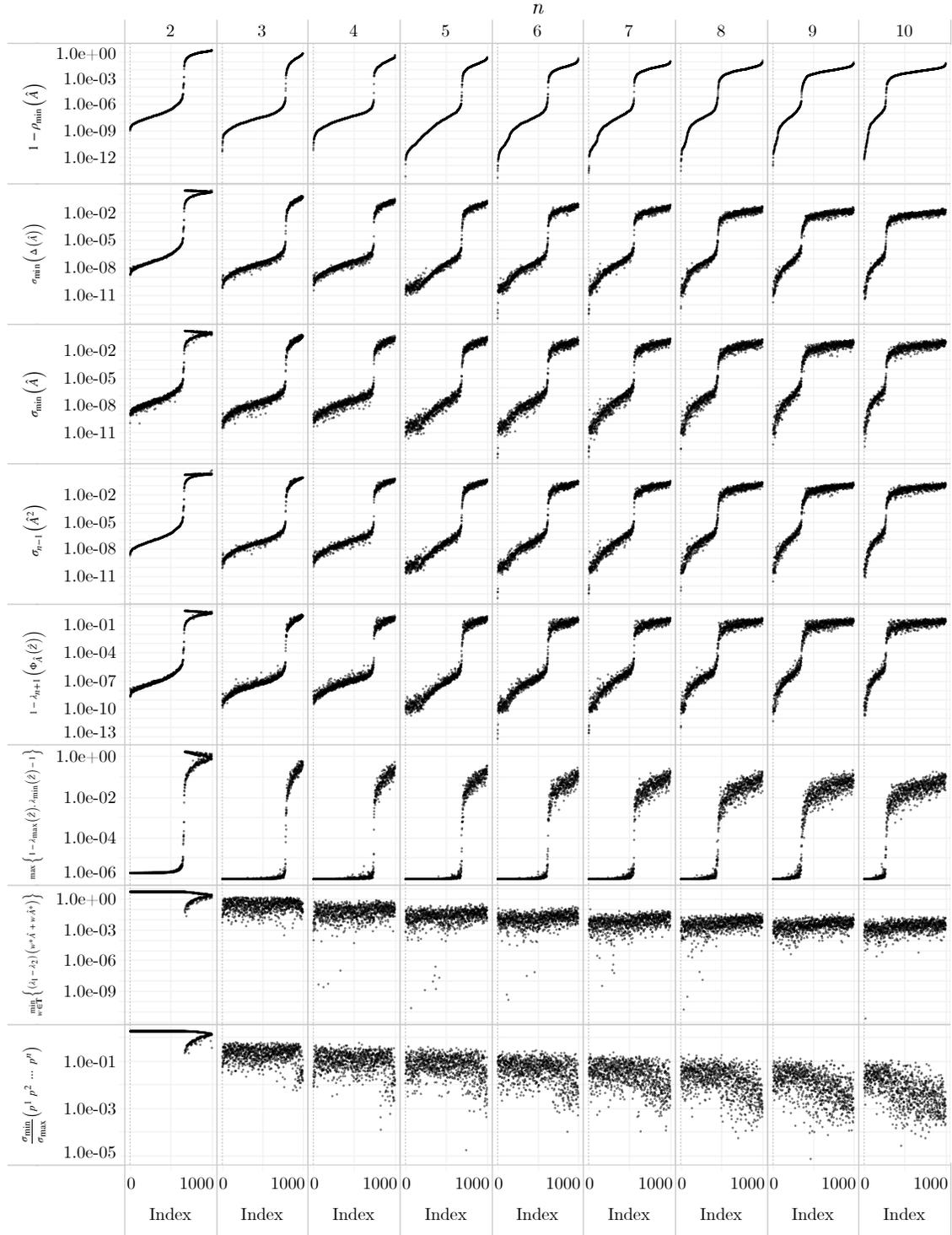}
\caption{Disk and non-disk matrices}
\label{plot}
\end{figure}

We consider the {\em proximal map} for the numerical radius,
$\mbox{prox}_r \colon \Mn \to \Mn$ defined by
\[
\mbox{prox}_r(Y) ~=~
\mbox{argmin} \big\{ r(A) + \lambda \|A-Y\|^2 : A \in \Mn \big\} \qquad (Y \in \Mn),
\]
where $\| \cdot \|$ denotes the Frobenius norm, and $\lambda > 0$ is a fixed constant.  For each dimension $n$, we generate a thousand independent random matrices $Y$ of unit norm by rescaling normally distributed matrices.  By solving a semidefinite program (described below), we then compute the matrix 
$A = \mbox{prox}_r(Y)$ and its numerical radius $r(A)$, along with several associated measures, each of which we plot in its own row in Figure \ref{plot}.

The first row measures how far the field of values of a matrix $A$ diverges from a disk.  By definition, we have
\[
\max_{z \in W(A)} \mbox{Re}(w^*z) ~\le~ r(A) \quad \mbox{for all}~ w \in \T
\]
where $\T$ denotes the unit circle in the complex plane $\C$, and disk matrices are characterized by equality for all $w$.  The left-hand side is the maximum eigenvalue of the Hermitian matrix $\frac{1}{2}(w^*A + wA^*)$, so calculating a natural measure of divergence,
\bmye \label{divergence}
1 ~-~ \frac{1}{r(A)} \min_{w \in \T} \max_{z \in W(A)} \mbox{Re}(w^*z),
\emye
is straightforward using the {\tt chebfun} package \cite{chebfun14}.  This measure vanishes for a nonzero matrix $A$ if and only $A$ is either a disk matrix, or a multiple of the identity (an exceptional case that we do not observe experimentally).

We sequence our thousand random matrices $A$ in ranked increasing order of the divergence measure and preserve this order in each of the subsequent rows of graphs, which plot this same ranking against the various measures for the scaled matrix $\hat A = \frac{1}{r(A)}A$.  The step functions in the first row of Figure \ref{plot}, involving several orders of magnitude, clearly suggest an underlying partition of the output set into disk and non-disk matrices.  The second row of Figure \ref{plot} plots another precise measure of divergence from the disk-matrix property, which we introduce shortly.  Numerical limitations inevitably blur the clarity of the step functions:  for comparison, for Figure \ref{plot2} in the Appendix, we repeat the calculations in dimensions $n=2,3$ to higher precision, refining the step-function behavior.

The experimental results also suggest that the random output $A$ is often singular, and indeed has a multiple zero eigenvalue.  We measure these latter properties, through singular values of $A$ and $A^2$, in the third and fourth rows of Figure \ref{plot} (which incidentally also verify that $A$ is never a multiple of the identity matrix).  In fact, in our experiments, these various properties, along with two others related to the associated semidefinite program and plotted in the fifth and sixth rows, appear almost surely equivalent.

We note that normalizing the random input $Y \in \Mn$  encourages the output 
$\mbox{prox}_r(Y)$ to be a disk matrix.  All sufficiently small $Y$ satisfy $\mbox{prox}_r(Y) = 0$:   specifically, this holds exactly when $\lambda Y$ lies in the unit ball for the norm dual to the radius.  On the other hand, since the set of disk matrices is a cone, the distance to this set for a scaled matrix $\alpha Y$ is proportional to the scaling $\alpha > 0$:  if $\alpha$ is large, 
$\mbox{prox}_r(\alpha Y)$ will not be a disk matrix.  The proportion of disk matrices seen depends on the balance between the normalization and the choice of the constant $\lambda$, which we fixed to be 
$\frac{3}{4}$ in the experiment.

\section{The semidefinite program}
Our calculation of the proximal map relies on the semidefinite representability of the numerical radius.  We denote the space of $n$-by-$n$ Hermitian matrices by $\Hn$, the cone of positive-semidefinite matrices by $\Hn_+$, and the corresponding semidefinite ordering by $\succeq$.  We consider the eigenvalues of a matrix $X \in \Hn$, listed by multiplicity, in nonincreasing order,
\[
\lambda_{\max}(X) = \lambda_1(X) \ge \lambda_2(X) \ge \cdots \ge \lambda_n(X) = \lambda_{\min}(X),
\]
and we use analogous notation for the singular values of an arbitrary matrix $Y$:
\[
\sigma_j(Y) = \sqrt{\lambda_j(Y^*Y)}.
\]
For a fixed matrix $A \in \Mn$, we define an affine map
$\Phi_A \colon \Hn \to {\mathbf H}^{2n}$ by
\[
\Phi_A(Z) ~=~ 
\left[
\begin{array}{cc}
Z &  A \\ 
A^* & - Z
\end{array}
\right]
\qquad
\mbox{for}~ Z \in \Hn.
\]
The following standard characterization of the numerical radius is due to Mathias \cite{mathias}, based in part on classical work of Ando \cite{ando}.  

\begin{thm}[SDP representation of the numerical radius] \label{mathias} 
\hspace{2cm} \mbox{} \hfill \mbox{}
The numerical radius of any matrix $A \in \Mn$ is the minimum value of the function 
$\lambda_{\max} \circ \Phi_A$, and the minimum is attained.
\end{thm}

\noindent
Any minimizer of the function $\lambda_{\max} \circ \Phi_A$ is called an {\em SDP-representor} for  $A$.

Using this result, we deduce that the minimum
\bmye \label{experiment}
\min_{A \in \Mn,~ Z \in \Hn}
\Big\{ 
\lambda_{\max} 
\left[
\begin{array}{cc}
Z &  A \\ 
A^* & - Z
\end{array}
\right]
+
\frac{1}{2} \|A-Y\|^2 
\Big\},
\emye
is attained, and at optimality, the matrix $A$ must be $\mbox{prox}_r(Y)$.  Both the functions $\lambda_{\max}$ and $\|\cdot\|^2$ are semidefinite representable, so we can easily convert this minimization problem into an equivalent semidefinite program:  in our numerical results, we solve it using the {\tt cvx} package.

We can check the disk matrix property for any matrix, and in particular for the output of our experiment, using straightforward numerical linear algebra, rather than relying on the divergence measure (\ref{divergence}) and the {\tt chebfun} package.  We apply the following result, which we derive later from work of Ando \cite{ando}. 

\begin{thm} \label{check-disk}
A matrix $A \in \Mn$ has field of values the unit disk if and only its numerical radius is one and there exists a nonzero polynomial $n$-vector $p(\cdot)$ of degree less than $n$ satisfying
\bmye \label{certificate-property}
(w^* A + w A^*) p(w) ~=~ 2p(w) \quad \mbox{for all}~ w \in \T.
\emye
\end{thm}

We call such a polynomial vector $p(\cdot)$ a {\em disk certificate}.  A disk certificate for a general nonzero matrix $A \in \Mn$ is just a disk certificate for $\frac{1}{r(A)} A$.  If the image $p(\T)$ spans the space $\Cn$, then we call the certificate {\em spanning}.

To convert this result into a numerical tool, we identify any polynomial $n$-vector $p(\cdot)$ of degree less than $n$ with its coefficient vectors:
\bmye \label{coefficients}
p(w) = \sum_{j=1}^n p^j w^{j-1} \qquad \mbox{for}~ p^j \in \Cn \quad (j=1,2,\ldots,n).
\emye
Certificates $p(\cdot)$ satisfying equation (\ref{certificate-property}) then correspond to nontrivial solutions of the recursion
\begin{eqnarray*}
Ap^{k+1} + A^* p^{k-1} &=& 2p^k \qquad (k=0,1,2,\ldots,n+1) \\
p^{-1} ~=~ p^0 &=& 0 ~=~ p^{n+1} ~=~ p^{n+2},
\end{eqnarray*}
which we can rewrite in matrix form,
\[
\Delta(A) 
\left[
\begin{array}{c}
p^1 \\ p^2 \\ p^3 \\ \vdots \\ p^n 
\end{array}
\right]
~=~
\left[
\begin{array}{c}
0 \\ 0 \\ 0 \\ \vdots \\ 0 
\end{array}
\right],
\]
for the $n(n+2)$-by-$n^2$ matrix
\bmye
\Delta(A) ~=~
\left[
\begin{array}{cccccc}
A & 0 & 0 & 0& \cdots & 0 \\
-2I & A & 0 & 0 & \cdots & 0 \\
A^* & -2I & A & 0 & \cdots & 0 \\
0 & A^* & -2I & A & \cdots & 0 \\
\vdots & \ddots & \ddots & \ddots & \ddots & \vdots \\
0 & 0 & \cdots & A^* & -2I & A \\
0 & 0 & \cdots & 0 & A^* & -2I \\
0 & 0 & \cdots & 0 & 0 & A^*
\end{array}
\right].
\emye
Thus, the disk property amounts simply to linear dependence of the columns of the matrix 
$\Delta(\frac{1}{r(A)}A)$ --- the second row of Figure \ref{plot} plots the smallest singular value of this matrix.  It is easy to see that a certificate $p(\cdot)$ is spanning exactly when its coefficient vectors $p^j$ comprise a basis for $\Cn$, a property that we can write in terms of a condition measure:
\bmye \label{condition}
\frac{\sigma_{\min}}{\sigma_{\max}} (p^1 ~ p^2 \cdots p^n) ~>~ 0.
\emye
The final row of Figure \ref{plot} plots this measure.

We call a matrix $A \in \Mn$ {\em separated} if the largest eigenvalue of the Hermitian matrix 
$w^*A + wA^*$ is simple for all complex units $w$, or in other words
\bmye \label{separation}
\min_{w \in \T} (\lambda_1 - \lambda_2)(w^*A + w A^*) ~>~ 0.
\emye
In that case, any disk certificate must be unique, up to nonzero scalar multiplication, or equivalently, the matrix on the left-hand side has a one-dimensional null space.  The penultimate row of Figure \ref{plot} plots this measure.  We call separated matrices with spanning certificates {\em strongly certified}.

As we already observed, the output of our experiment, $A = \mbox{prox}_r(Y)$ is often a disk matrix.  In those instances, we can determine whether $A$ is strongly certified by checking the conditions (\ref{condition}) and (\ref{separation}).  Indeed, we find, numerically, that the disk matrices generated are always strongly certified.

\section{Disk certificates} 
We begin our more formal development by justifying our definition of a disk certificate.  For any matrix $A \in \Mn$ and complex unit $w \in \T$, we have
\[
2\max \{ \mbox{Re}(w^*z) : z \in W(A) \}
~=~
\lambda_{\max} (w^* A + w A^*).
\]
Hence $A$ is a disk matrix if and only if, for all $w \in \T$, there exists a nonzero vector 
$p(w) \in \Cn$ satisfying 
\bmye \label{disk-certificate}
(w^* A + w A^*) p(w) ~=~ 2r(A)p(w).
\emye
However, using two key results, due to Ando and Crouzeix, which we summarize next, we can restrict attention to vectors $p(w)$ satisfying equation (\ref{disk-certificate}) whose components are polynomials in $w$ of degree less than $n$.  To see this, we rely on the following restatement of \cite[Theorem 1]{ando}.

\begin{thm} \label{ando}
A matrix $A \in \Mn$ has numerical radius $r(A) \le 1$ if and only if it has an {\bf Ando representation}
$A = 2SUC$, for matrices $S,C \in \Hn_+$ satisfying $S^2 + C^2 = I$, and a unitary matrix $U$.
\end{thm}

\noindent
We deduce the following result, subsuming Theorem \ref{check-disk}.

\begin{thm} \label{certificate}
Suppose that a matrix $A \in \Mn$ with unit numerical radius has Ando representation $A=2SUC$, as in Theorem \ref{ando}.  Then disk certificates for $A$ are exactly nonzero polynomial $n$-vectors $p$ of degree less than $n$ satisfying
\bmye \label{null}
(UC-wS)p(w) = 0 ~~\mbox{for all}~ w \in \C.
\emye
Furthermore, the following properties are equivalent:
\begin{enumerate}
\item[{\rm (i)}]
The matrix $UC-wS$ is singular for all $w \in \C$;
\item[{\rm (ii)}]
$A$ has a disk certificate;
\item[{\rm (iii)}]
$A$ is a disk matrix. 
\end{enumerate}
\end{thm}

\pf
The first claim follows directly from \cite[Proposition 8.2]{lewis-overton-disks}.
The equivalence of properties (i) and (iii) follows \cite[Lemma 10.1]{Cro16}.  Clearly the existence of a disk certificate implies the singularity property (i).  Conversely, if the singularity property (i) holds then the matrix pencil $UC-wS$ is singular.  The Kronecker form for singular matrix pencils \cite[pp.~106--108]{vandooren} now guarantees the existence of a nonzero polynomial $n$-vector $p$ of degree less than $n$ satisfying equation (\ref{null}), and as we have already seen, this $p$ is a disk certificate.
\finpf

\noindent
By rescaling, we deduce the justification for our disk certificate definition.

\begin{cor}
A matrix is a disk matrix if and only if it has a disk certificate.
\end{cor}

To consider more specifically the set of strongly certified matrices, we begin with the following easy result.

\begin{prop}
In the space $\Mn$, the set of separated matrices is open.
\end{prop}

\pf
Consider a sequence of non-separated matrices $A_k \in \Mn$, for $k=1,2,3,\ldots$, with limit $A \in \Mn$.  By assumption, there exist corresponding complex units $w_k \in \C$ such that the largest two eigenvalues of the matrix $w_k^* A_k + w_k A_k^*$ are equal.  Taking a subsequence, we can suppose that $w_k$ converges to some complex unit $w \in \C$.  By continuity, the largest two eigenvalues of the matrix 
$w^* A + wA^*$ are also equal, so the result follows.
\finpf

\begin{thm} \label{relative}
Within the set of all disk matrices in $\Mn$, the set of strongly certified disk matrices is relatively open.
\end{thm}

\pf
Consider a strongly certified disk matrix $A \in \Mn$.  We claim that all nearby disk matrices are strongly certified.  If not, then there exists a sequence of disk matrices $A_k \to A$, each of which is not not strongly certified.  By the previous result, each $A_k$ is separated, so its disk certificate $p_k$ (unique up to scalar multiplication) is not spanning.  By scaling, we can suppose that each certificate has unit norm.  (Here, we define the norm of any polynomial vector $p$, as in equation (\ref{coefficients}), as the Frobenius norm of the matrix with columns its coefficients $p^j$.)  Taking a subsequence, we can suppose that $p_k$ converges to some nonzero limit polynomial $p$, which by continuity must be the unique certificate for $A$, and furthermore is not spanning.  This contradiction completes the proof.
\finpf

\begin{exa}[A strongly certified disk matrix] \label{two}
{\rm
The matrix
\[
A ~=~
\left[
\begin{array}{cc}
0 & 2 \\
0 & 0 
\end{array}
\right]
\]
has field of values the unit disk, since it has the Ando representation
\[
2
\left[
\begin{array}{ccc}
1 & 0 \\
0 & 0 
\end{array}
\right]
\left[
\begin{array}{ccc}
0 & 1 \\
1 & 0
\end{array}
\right]
\left[
\begin{array}{ccc}
0 & 0 \\
0 & 1
\end{array}
\right],
\]
and the polynomial vector 
\[
\left[
\begin{array}{c}
1 \\ w
\end{array}
\right]
\]
is a spanning disk certificate.  Furthermore, the matrix $A$ is separated, since the matrix $w^*A + wA^*$ has the distinct eigenvalues $\pm 2$, so $A$ is strongly certified.
}
\end{exa}

\begin{exa}[A disk matrix with no spanning certificate] \label{three}
{\rm
The matrix
\[
A ~=~
\left[
\begin{array}{ccc}
0 & 0 & 2 \\
0 & \frac{4}{5} & 0 \\
0 & 0 & 0 
\end{array}
\right].
\]
has field of values the unit disk, since it has the Ando representation
\[
2
\left[
\begin{array}{ccc}
1 & 0 & 0 \\
0 & \frac{1}{\sqrt{2}} & 0 \\
0 & 0 & 0
\end{array}
\right]
\left[
\begin{array}{ccc}
0 & 0 & 1 \\
\frac{3}{5} & \frac{4}{5} & 0 \\
\frac{4}{5} & -\frac{3}{5} & 0
\end{array}
\right]
\left[
\begin{array}{ccc}
0 & 0 & 0 \\
0 & \frac{1}{\sqrt{2}} & 0 \\
0 & 0 & 1
\end{array}
\right],
\]
and the polynomial vector
\[
\left[
\begin{array}{c}
1 \\ 0 \\ w
\end{array}
\right]
\]
is a disk certificate.  However, a direct application of Theorem \ref{certificate} shows that 
any disk certificate for $A$ must have second component zero, and hence $A$ has no spanning certificate. 
}
\end{exa}

Singular matrices need not be disk matrices, and as Example \ref{three} highlights, disk matrices need not be strongly certified.  However, no such examples emerge from our computational experiment:  whenever the output matrix $A$ is singular, it is a strongly certified disk matrix.  The numerical radius must be partly smooth around such matrices, as a consequence of the following result, which, in the light of Theorem \ref{relative}, is just a restatement of \cite[Theorem 7.4]{lewis-overton-disks}.

\begin{thm}
The set of all strongly certified disk matrices in $\Mn$ is a real-analytic manifold of real codimension $2n$, with respect to which the numerical radius is partly smooth.
\end{thm}

\section{Random disk matrices and SDP-representors}
For any matrix $A \in \Mn$, Theorem \ref{mathias} ensures the existence of an SDP-representor, or in other words a minimizer of the function 
$\lambda_{\max} \circ \Phi_A$.  Our computational experiment (\ref{experiment}) suggests a variety of interesting properties of SDP-representors.  Each random instance solves a single semidefinite program which, by decomposing into two subproblems, we can view as an experiment on random disk matrices:
\begin{itemize}
\item
Compute $A = \mbox{prox}_r(Y)$ for a random input matrix $Y \in \Mn$;
\item
Minimize $\lambda_{\max} \circ \Phi_A$ to find an SDP-representor $Z$ for $A$.
\end{itemize}
We note that the SDP-representor $Z$ may not be unique.

In the experiment, the following properties often hold.  Furthermore, in each random instance, with occasional numerical anomalies, either all the properties hold, or none do.
\begin{enumerate}
\item[(i)]
$A$ is singular.
\item[(ii)]
$A^2$ has rank strictly less than $n-1$.
\item[(iii)]
$\lambda_{\max}(Z) = r(A)$ and $\lambda_{\min}(Z) = -r(A)$.
\item[(iv)]
$A$ is a disk matrix.
\item[(v)]
The largest eigenvalue of $\Phi_A(Z)$ has multiplicity strictly larger than $n$.
\item[(vi)]
$A$ is a strongly certified disk matrix.
\end{enumerate}
In light of property (v), we say that a matrix $A$ has the {\em multiplicity property} if there exists an SDP-representor $Z$ (not necessarily unique) such that the largest eigenvalue of $\Phi_A(Z)$ has multiplicity strictly larger than $n$.

We can quantify each of the first five properties with a corresponding nonnegative measure that vanishes exactly when the property holds:
\begin{itemize}
\item
$\sigma_{\min}(A)$
\item
$\sigma_{n-1}(A^2)$
\item
$\max\big\{ r(A) - \lambda_{\max}(Z) , \lambda_{\min}(Z) - r(A) \big\}$
\item
$\sigma_{\min}\big(\Delta\big(\frac{1}{r(A)}A\big)\big)$
\item
$r(A) - \lambda_{n+1}\big(\Phi_A(Z)\big)$.
\end{itemize}
We verify the strong certification (property (vi)) by checking the conditions (\ref{condition}) and (\ref{separation}) in conjunction with property (iv).  To illustrate the coincidence of the six properties, for each matrix size $n=1,2,3,\ldots,$ we order the random instances in increasing order of the disk-property divergence measure (\ref{divergence}), and then simultaneously plot the measures, scaled appropriately, in Figure \ref{plot}.

Of course, for {\em arbitrary} matrices $A$, the six properties are far from equivalent.
We shall see (in Corollary \ref{corollary}) the general implications
\bmye \label{implications1}
\mbox{(iv)} ~\Rightarrow~ \mbox{(iii)} ~\Rightarrow~ \mbox{(ii)} ~\Rightarrow~ \mbox{(i)}. 
\emye
Furthermore, when the SDP-representor $Z$ is unique, we shall see (in Theorems \ref{sufficient} and \ref{necessary}) the implications
\[
\mbox{(vi)} ~\Rightarrow~ \mbox{(v)} ~\Rightarrow~ \mbox{(iv)}. 
\]
On the other hand, examples show the relationships
\[
\mbox{(i)} ~\not\Rightarrow~ \mbox{(ii)} ~\not\Rightarrow~ \mbox{(iii)} ~\not\Rightarrow~ \mbox{(iv)}
~\not\Rightarrow~ \mbox{(v)}.
\]
The diagonal matrix $A=\mbox{Diag}(0,1)$ illustrates the first relationship, and a direct calculation show that the diagonal matrices $A=\mbox{Diag}(0,0,1)$ and $Z=\mbox{Diag}(1,1,0)$ illustrate the second.  For the third relationship, we present an example after Corollary \ref{corollary}, and Example \ref{many} illustrates the fourth.  The precise relationship between the multiplicity property and strong certification remains unclear.

To explore these properties, we begin by characterizing SDP-representors.  We first consider matrices with unit numerical radius:  the general case follows easily by scaling. 

\begin{thm} \label{horn}
A matrix $Z \in \Hn$ is an SDP-representor for a matrix $A \in \Mn$ with unit numerical radius if and only if 
\bmye \label{contraction}
I_n \succeq Z \succeq -I_n \quad \mbox{and} \quad
A ~=~ (I-Z)^{\frac{1}{2}} U (I+Z)^{\frac{1}{2}} ~ \mbox{for a contraction}~ U \in \Mn.
\emye
Furthermore, there exist such matrices $Z$ and $U$, with $U$ in fact unitary.
\end{thm}

\pf
By Theorem \ref{mathias}, a matrix $Z$ is an SDP-representor for $A$ if and only if the matrix
$I_{2n} - \Phi_A(Z)$ is positive semidefinite, and this is equivalent to the property (\ref{contraction}) by \cite[Theorem 7.7.9]{HorJoh13}.  The final result is proved in \cite{ando}.
\finpf

The proof of the following result essentially follows \cite[Lemma 10.1]{Cro16}.

\begin{thm} \label{lemma}
If a matrix $A \in \Mn$ with field of values the unit disk satisfies the property (\ref{contraction}), then the matrix
\[
U (I+Z)^{\frac{1}{2}} - w(I-Z)^{\frac{1}{2}}
\]
is singular for all $w \in \C$.
\end{thm}

\pf
Consider the two positive semidefinite Hermitian matrices $S = (I-Z)^{\frac{1}{2}}$ and 
$C = (I+Z)^{\frac{1}{2}}$.  By assumption, for any $w \in \C$ satisfying $|w|=1$ there exists a unit vector $v \in \Cn$ satisfying $v^*Av = w$.  We now observe the inequalities
\[
1 = |v^*Av| = 2|v^*SUCv| \le 2 \|Sv\| \|UCv\| \le 2 \|Sv\| \|Cv\| \le \|Sv\|^2 + \|Cv\|^2 = 1,
\]
Equality in the first inequality ensures $UCv = w' Sv$ for some $w' \in \C$.  Equality in the last inequality, along with the final equality, ensures $\|Sv\| = \frac{1}{\sqrt{2}} = \|Cv\|$.  We deduce
\[
w = v^*Av = 2v^*SUCv = 2v^*Sw' Sv = w',
\]
so $UCv = w Sv$, and hence the matrix $UC-wS$ is singular.  The result follows.
\finpf

\noindent
We deduce the implications (\ref{implications1}), which we summarize in the following result.

\begin{cor} \label{corollary}
If a matrix $A \in \Mn$ has an SDP-representor with both the eigenvalues $\pm r(A)$, as holds in particular if $A$ is a disk matrix, then $\mbox{rank}(A^2) < n-1$, and hence $A$ is singular.
\end{cor}

\pf
By scaling, we can assume $r(A)=1$.  By Theorem \ref{horn}, for any SDP-representor $Z$, there exists a contraction $U \in \Mn$ satisfying the property (\ref{contraction}).  We deduce
\[
A^2 ~=~ (I-Z)^{\frac{1}{2}} U (I-Z^2)^{\frac{1}{2}} U (I+Z)^{\frac{1}{2}},
\]
and the rank condition now follows from the fact that $\mbox{rank}(I-Z^2) < n-1$.  In the case when $A$ is a disk matrix, we apply Theorem \ref{lemma}, and consider the cases $w=0$ and $w \to \infty$.
\finpf

The converse of Corollary \ref{corollary} is false.  For example, it is easy to verify that if we define
\[
A ~=~
\left[
\begin{array}{ccc}
0 & 0 & 1 \\
0 & 1 & 0 \\
0 & 0 & 0
\end{array}
\right]
\qquad
\mbox{and}
\qquad
Z ~=~
\left[
\begin{array}{ccc}
1 & 0 & 0 \\
0 & 0 & 0 \\
0 & 0 & -1
\end{array}
\right],
\]
then $A$ has unit numerical radius, $Z$ is an SDP-representor with eigenvalues $\pm 1$, and yet $A$ is not a disk matrix.

\begin{prop} \label{optimal}
If matrix $A \in \Mn$ with unit numerical radius has Ando representation $A=2SUC$, as in Theorem \ref{ando}, then the matrix
\[
Z ~=~ I - 2S^2 ~=~ 2C^2 - I
\]
is an SDP-representor, and the largest eigenvalue of the matrix $\Phi_A(Z)$
has multiplicity at least $n$.
\end{prop}

\pf
The matrix
\bmye \label{factorization}
B ~=~
I-\Phi_A(Z) ~=~
2
\left[
\begin{array}{c}
S \\
-C U^*
\end{array}
\right]
\begin{array}{c}
\left[ S ~~ -UC \right] \\
\mbox{}
\end{array}
\emye
is positive semidefinite.  We deduce the inequality
$\lambda_{\max}\big(\Phi_A(Z)\big) \le 1$, 
and in fact equality must hold, by Theorem \ref{mathias}, and $Z$ must be an SDP-representor.  Furthermore, the multiplicity of the largest eigenvalue of $\Phi_A(Z)$ is $2n - \mbox{rank}\,B$.  Since the rank of any matrix $Y$ coincides with the rank of the matrix $Y^* Y$, we see $\mbox{rank}\, B \le n$, and the result now follows.
\finpf

\noindent
Turning to the multiplicity property, we first note that the following corollary of the previous result.

\begin{cor}
Any matrix $A \in \Mn$ has an SDP-representor $Z \in \Hn$ for which the largest eigenvalue of $\Phi_A(Z)$ has multiplicity at least $n$.
\end{cor}

Consider, for example, Example \ref{two}.  Any SDP-representor $Z \in {\mathbf H}^2$ has both the eigenvalues $\pm 1$, by Corollary \ref{corollary}, and hence $Z^2 = I$.  By Theorem \ref{horn}, there exists a contraction 
$U \in {\mathbf M}^2$ such that 
\[
A ~=~ (I-Z)^{\frac{1}{2}}U(I+Z)^{\frac{1}{2}}, 
\]
so
\[
A(I-Z)~=~ (I-Z)^{\frac{1}{2}}U(I-Z^2)^{\frac{1}{2}}(I-Z)^{\frac{1}{2}} ~=~ 0.
\]
Hence the matrix $I-Z$ must be diagonal, whence so is $Z$.
We deduce that the unique SDP-representor is the matrix 
\[
Z ~=~ 
\left[
\begin{array}{cc}
-1 & 0 \\
0 & 1 
\end{array}
\right].
\]
In that case, the matrix
\[
\Phi_A(Z)
~=~
\left[
\begin{array}{cccc}
-1 & 0 & 0 &  2 \\
 0 & 1 & 0 &  0 \\
 0 & 0 & 1 &  0 \\
 2 & 0 & 0 & -1
\end{array}
\right]
\]
has eigenvalues (by multiplicity) $1,~1,~1$, and $-3$, illustrating property (v).

\begin{que}
{\rm
Must SDP-representors of strongly certified disk matrices be unique?
}
\end{que}

\section{Unique SDP-representors}
As we noted in the previous section, the matrices $A$ generated by our random experiment are disk matrices exactly when they satisfy the multiplicity property.  In this section we explain this phenomenon, at least in the case when the SDP-representor is unique.  We begin with the following result.

\begin{thm}[Sufficient condition for disk matrices] \label{sufficient}~~
Consider a matrix \mbox{$A \in \Mn$} with a unique SDP-representor $Z \in \Hn$.  If $A$ has the multiplicity property, then it must be a disk matrix.
\end{thm}

\pf
After scaling, we may assume $r(A) = 1$.  Given an Ando representation $A=2SUC$, as in Theorem \ref{ando}, Proposition \ref{optimal} implies $Z = I - 2S^2$.
Using the multiplicity property, and following the argument of Proposition \ref{optimal}, we deduce, for all numbers $w \in \C$,
\begin{eqnarray*}
\lefteqn{n ~>~ \mbox{rank}\big(I-\Phi_A(Z)\big) ~=~ \mbox{rank} [ S ~~ -UC ]} \\
& & ~\ge~ 
\mbox{rank} 
\left(
\begin{array}{c}
[ S ~~ -UC ] \\
\mbox{}
\end{array}
\left[
\begin{array}{c}
w I \\
I
\end{array}
\right]
\right)
~=~
\mbox{rank}(UC-wS),
\end{eqnarray*}
so the matrix $UC-wS$ is singular.  The result now follows by Theorem \ref{ando}.
\finpf

\begin{que}
{\rm
Is the uniqueness assumption necessary?  In other words, for a nondisk matrix $A$ with an SDP-representor $Z$, can the largest eigenvalue of the matrix $\Phi_A(Z)$ have multiplicity exceeding $n$?
}
\end{que}

\begin{que} \label{converse}
{\rm
Is the converse true? That is, if a disk matrix $A$ has a unique SDP-representor $Z$, can the largest eigenvalue of the matrix $\Phi_A(Z)$ have multiplicity $n$?
}
\end{que}

Routine algebra, and the classification in \cite{lewis-overton-disks} verifies the converse of Theorem \ref{sufficient} in the $2$-by-$2$ and $3$-by-$3$ cases.  On the other hand, it can fail without the uniqueness assumption, as shown by the following example.

\begin{exa}[A disk matrix with many SDP-representors] \label{many} \hfill \mbox{} \\
{\rm 
Consider the disk matrix $A$ of Example \ref{three}, which has no strong (or indeed spanning) certificate.  For any real number $s$, consider the matrix
\[
Z ~=~
\left[
\begin{array}{ccc}
-1 & 0 & 0 \\
0 & s & 0 \\
0 & 0 & 1 
\end{array}
\right].
\]
The matrix 
\[
\Phi_A(Z)
~=~
\left[
\begin{array}{cccccc}
-1 & 0 & 0 & 0 & 0 & 2 \\
0 & s & 0 & 0 & \frac{4}{5} & 0 \\
0 & 0 & 1 & 0 & 0 & 0 \\
0 & 0 & 0 & 1 & 0 & 0 \\
0 & \frac{4}{5} & 0 & 0 & -s & 0 \\
2 & 0 & 0 & 0 & 0 & -1 
\end{array}
\right]
\]
is unitarily similar (by row and column permutations) to 
\[
\left[
\begin{array}{cc}
1 & 0 \\
0 & 1
\end{array}
\right]
\oplus
\left[
\begin{array}{cc}
-1 & 2 \\
2 & -1
\end{array}
\right]
\oplus
\left[
\begin{array}{cc}
s & \frac{4}{5} \\
\frac{4}{5} & -s
\end{array}
\right],
\]
and hence has eigenvalues
\[
1,~ 1,~ 1, -3,~ \pm\sqrt{s^2 + \Big(\frac{4}{5}\Big)^2}.
\]
Thus $Z$ is an SDP-representor providing $|s| \le \frac{3}{5}$, and if $|s| < \frac{3}{5}$, then the multiplicity of the largest eigenvalue of $\Phi_A(Z)$ is three.
}
\end{exa}

On the other hand, in this example, setting $s = \pm \frac{3}{5}$ does produce an SDP-representor $Z$ for which the multiplicity (four) of the largest eigenvalue of the matrix $\Phi_A(Z)$ exceeds $n=3$.  This suggests the following question, more general than Question \ref{converse}.

\begin{que}
{\rm
Do all disk matrices $A$ satisfy the multiplicity property?
}
\end{que}

\noindent
Again, routine algebra and \cite{lewis-overton-disks} answers this question affimatively in the $2$-by-$2$ and $3$-by-$3$ cases.

Unlike the disk matrix $A$ of Examples \ref{three} and \ref{many}, the disk matrices produced by our random experiment are always strongly certified.  In that case we have a more positive result, as we argue next.  We begin with an elementary tool.

\begin{lem} \label{first}
The components of polynomial $k$-vector $q(w)$ (in the variable $w \in \C$) are linearly independent over $\C$ if and only its range $q(\C)$ spans the complex vector space $\C^k$.
\end{lem} 

\pf
The conclusion fails if and only if there exists a nonzero vector $a \in \C^k$ satisfying 
$\sum_j a_j q_j(w) = 0$ for all $w \in \C$, so the result follows.
\finpf

\noindent
This result extends, as follows.

\begin{lem} \label{second}
The dimension in the complex vector space $\C^k$ of the range of a polynomial $k$-vector $q(w)$ (in the variable $w \in \C$) equals the dimension of the span over $\C$ of its component polynomials.
\end{lem} 

\pf
Denote by $d$ the second dimension.  Denote by $q^d$ and $q^{k-d}$ the polynomial vectors obtained from $q$ by deleting the last $k-d$ and first $d$ components respectively.  Without loss of generality, we can suppose that $q^d$ has linearly independent components and that $q^{k-d} = Eq^d$ for some complex matrix $E$.  By Lemma \ref{first}, the range $q^d(\C)$ has dimension $d$ in the space $\C^d$, and the range $q(\C)$ is its image under the injective linear map $\binom{I}{E}$, so has the same dimension.
\finpf

\begin{thm}[Necessary condition for disk matrices]~ \label{necessary}
If a nonzero matrix \mbox{$A \in \Mn$} has a spanning disk certificate (or in particular is strongly certified), then it satisfies the multiplicity property.
\end{thm}

\pf
After scaling, we can assume that the matrix $A$ has unit numerical radius and an Ando representation $A=2SUC$, as in Theorem \ref{ando}, with a spanning disk certificate $p$.  By Proposition \ref{optimal}, the matrix $Z=I-2S^2$ is an SDP-representor, and by equation (\ref{factorization}) the vector
\[
q(w) ~=~  
\left[
\begin{array}{c}
wp(w) \\ p(w)
\end{array}
\right]
\]
is an eigenvector corresponding to the largest eigenvalue (namely one) of the matrix $\Phi_A(Z)$ for all $w \in \C$.  By assumption, the $n$ components of the vector polynomial $p$ are linearly independent.  The vector polynomial $q$ must therefore have a linearly independent set of more than $n$ components, since at last one component has degree strictly larger that any component of $p$.  By Lemma~\ref{second}, the matrix $\Phi_A(Z)$ therefore has a linearly independent set of more than $n$ eigenvectors corresponding to the eigenvalue $1$, so this eigenvalue has multiplicity larger than $n$.
\finpf 

\noindent
To conclude, in the case of a unique SDP-representor $Z$, Theorems \ref{sufficient} and \ref{necessary} together indicate a strong correlation between the disk property for matrices $A$ and the multiplicity property, exactly as we noted empirically.

\bibliographystyle{plain}
\small
\parsep 0pt

\def\cprime{$'$} \def\cprime{$'$}

\newpage
\section{Appendix}
\begin{figure}[H] 
\centering
\includegraphics[width=.9\textwidth]{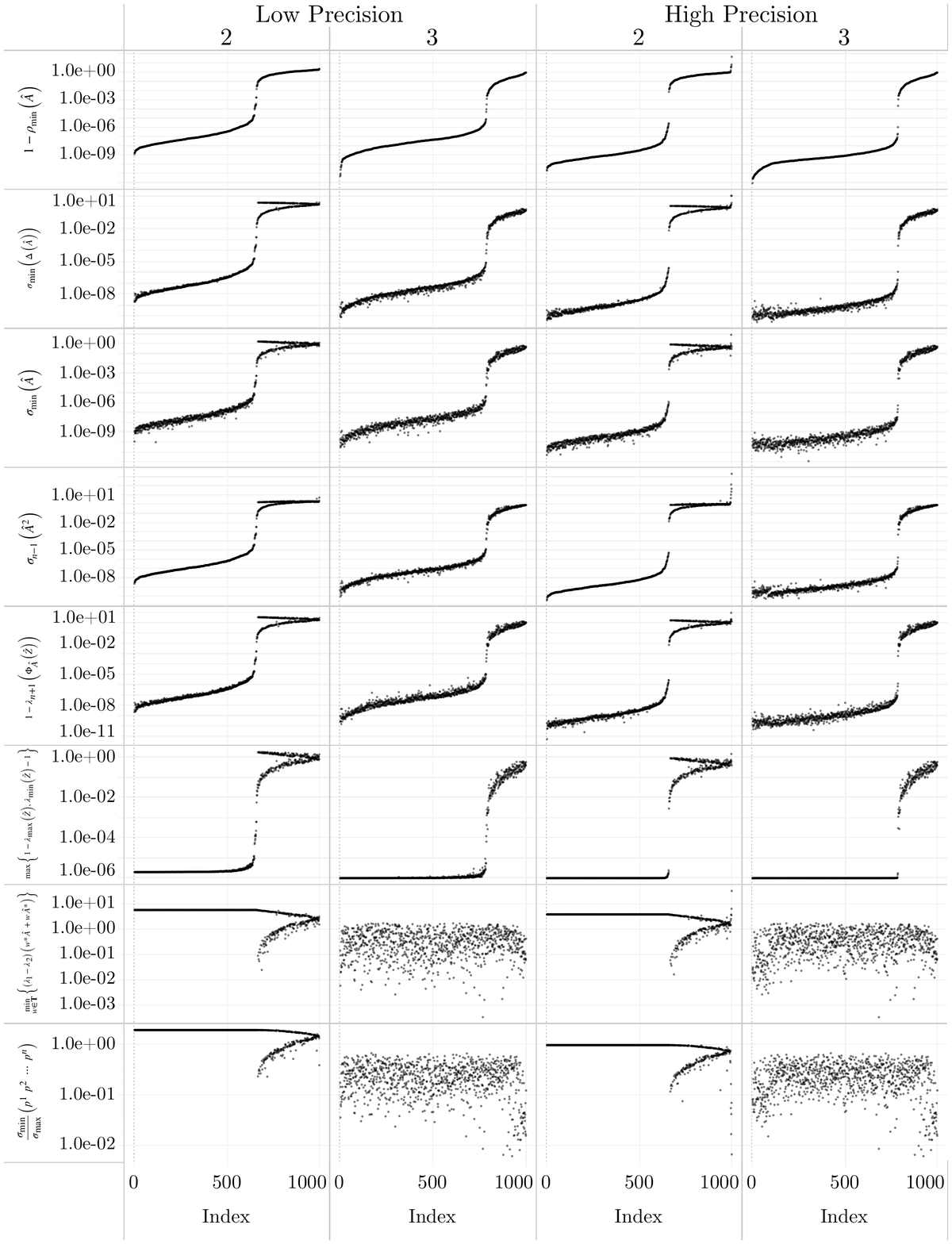}
\caption{Disk and non-disk matrices in higher precision}
\label{plot2}
\end{figure}

\end{document}